\documentclass[a4paper,10pt]{article}
\pdfoutput=1

\usepackage{amsmath, amsfonts, amssymb, amsthm}

\usepackage{graphicx}

\usepackage{xspace}
\usepackage[margin=1.0in]{geometry}
\usepackage{setspace}
\newtheorem{thm}{Theorem}
\newtheorem{lemma}{Lemma}

\author{Lucas Blakeslee}

\title{On the Genus of Random Regular Graphs}

\begin{document}

\begin{center}
\fontsize{13pt}{10pt}\selectfont
    \textsc{\textbf{ON THE GENUS OF RANDOM-REGULAR GRAPHS}}
    \end{center}
\vspace{0.1cm}
\begin{center}
   \fontsize{12pt}{10pt}\selectfont
    \textsc{Lucas Blakeslee}
\end{center}
\vspace{0.2cm}

\abstract{The genus of a graph is a topological invariant that measures the minimum genus of a surface on which the graph can be embedded without any edges crossing. Graph genus plays a fundamental role in topological graph theory, used to classify and study different types of graphs and their properties. We show that, for any integer $d \geq 2$, the genus of a random $d$-regular graph on $n$ nodes is $\frac{(d - 2)}{4}n(1 - \varepsilon) $ with high probability for any $\varepsilon > 0$}.

\section{Introduction}

Topological graph theory is a branch of mathematics that studies the topological properties of graphs and their embeddings on surfaces. It combines techniques from topology, algebraic geometry, and combinatorics to study the structure and properties of graphs. This field of mathematics is particularly useful for understanding the global properties of graphs, such as connectivity, and for studying the way in which graphs can be embedded in the plane or in more complicated topological spaces.

One central concept in topological graph theory is that of graph genus. The classification theorem for compact surfaces, first proposed by Möbius in the early 1860s \cite{Mobius1861} (though not proved rigorously until the 1920s by Brahana \cite{Brahana1921}), states that all closed surfaces are homeomorphic to spheres with some number of handles or crosscaps attached. We denote a sphere with $g$ handles attached $S_g$. This $g$ is called the \emph{genus} of the surface. The orientable genus of a graph $G$, $\gamma(G)$, is a topological invariant denoting the smallest number $g$ such that $G$ can be drawn onto the orientable surface $S_g$ without any edge crossings. In the simplest case, a graph with genus $0$ --- i.e., a planar graph --- can be drawn on the plane without edge crossings. A graph of genus 1 can be drawn without edge crossings on a torus, a graph of genus 2 on a surface with two handles, and so forth.

Much of the literature in topological graph theory focuses on the study of the embeddings of certain families of graphs, such as complete graphs \cite{Ringel1969}, wheel graphs \cite{Chen2018}, and complete bipartite (and more generally, $k$-partite) graphs \cite{Bouchet1978, Bettayeb2010}, to name a few. Many graph problems, such as the problem of finding a maximal independent set, finding a minimum $t$-spanner, the sparsest-cut problem, and several constraint satisfaction problems, are NP-hard for graphs of genus greater than zero \cite{Baker1983, Abboud2020, Kobayashi2018, Cai2017}. Therefore, the genus of a graph can be an important measure of complexity in analyzing the efficiency of graph algorithms. Determining the genus of an arbitrary graph is hard, however. In fact, we know the graph genus problem to be NP-complete \cite{Thomassen1989}.

The goal of this paper is to characterize the growth of the typical genus of a random regular graph, which can range from 0 in the case of planar graphs to the genus of the complete graph $K_n$, $\left\lceil \frac{(n - 3)(n-4)}{12} \right\rceil$, as per \cite{Ringel1969}. Related work has been performed by Asplund et al. \cite{Asplund2017}, who investigate the so-called $k$-planar crossing number of random regular graphs, another metric of how ``distant'' a graph is from being planar; however, the $k$-planar crossing number of a graph and its genus can be arbitrarily far apart. Additionally, Archdeacon and Grable have showed that the genus of the Erdős-Rényi graph grows as $\frac{n^2}{24}$ \cite{Archdeacon1995}. However, as we show here, the genus of a random $d$-regular graph grows with high probability as $\frac{d - 2}{4}n$ where $n$ is the number of nodes.

\section{Results}

We begin by bounding the number of cycles on the graph, which provides a bound on the number of faces, which in turn can then be used in conjunction with Euler's formula to describe the asymptotic growth of the genus. We say that an event occurs with high probability (w.h.p.) if it's probability goes to 1 as $n \to \infty$. For a function $f(G)$ of a random $d$-regular graph and a function $h(n)$, we say that $f(G) = o(h(n))$ w.h.p. if $f(G) < ch(n)$ w.h.p. for any $c > 0$. We will use $h(n) = n$ and $h(n) = 1$.

\begin{lemma}
\label{thm:one}
In a random $d$-regular graph, for any positive constant $m$, there are fewer than $cn$ cycles of length $m$ or less w.h.p. for any constant $c > 0$.
\end{lemma}

\emph{Proof.} We will first show that the expected number of such cycles is $O(1)$.

Markov's inequality states that if $X$ is a non-negative random variable, then for all $t$,
\begin{center}
    ${Pr}\left[X \geq t \right] \leq \mathbb{E}\left[ X \right] / t$.
\end{center}

\noindent By Markov's inequality, the probability that the number of said cycles exceeds $cn$ for any constant $c > 0$ is thus at most $O(\frac{1}{n}) = o(1)$.

There are $\frac{t!}{2t}$ possible ways to connect $t$ vertices in a cycle. For sufficiently large $n$, due to the fact that $m$ is a constant, even if we condition upon there being $m$ edges, there are still almost $dn$ ``half edges'' available, meaning the probability that a given set of t vertices forms a cycle in a particular order is less than $(d/n)^t$.

We sum over all possible orderings, and arrive on a bound for the number of cycles of length $\leq m$ of

\vspace{-8mm}

\begin{center}
    \[ \sum_{t=2}^m {n \choose t} \frac{t!}{2t} (d/n)^{t} \]
\end{center}

\noindent This expression is $O(1)$ if the term $m = O(1)$ (in fact, it is at most $m$ times its largest term, $\frac{d^m}{2m}$), from which the proof. \hfill $\square$






\begin{thm}
\label{thm:two}
For $d \geq 2$, the genus of a random $d$-regular graph on $n$ nodes grows w.h.p. as $\frac{(d - 2) n}{4}$
\end{thm}

\emph{Proof.} For a given random-regular graph, let $\alpha$ denote $g/n$, where $g$ is the graph's genus and $n$ is the number of nodes.

Euler's formula states that:

\begin{center}
    $V - E + F = 2 - 2g$
\end{center}

\noindent where $V$ is the number of vertices, $E$ the number of edges, and $F$ the number of faces in the embedding. Since $V = n$ and $E = dn/2$,

\begin{align*}
    &1 - d/2 + F/n = 2/n - 2\alpha \\
    \vspace{3mm}
    \implies \alpha &= \frac{d - 2}{4} + \frac{1}{n} - \frac{F}{2n} \\
    \vspace{3mm}
    &= \frac{d - 2}{4} + \frac{1}{n} - \left( \frac{d}{4} \right) \left( \frac{F}{E} \right) \tag{$\star$} \label{eq:one} \\
\end{align*}

To show our result, we will show that $F/E$ is $o(1)$ w.h.p. Let $c(m)$ denote the number of faces of length $m$ or less. Since each of these faces includes at least one edge, the other $F - c(m)$ faces have at least $m$ edges each, and each edge can appear in at most two faces, we have:

\begin{center}
    $2E \geq c(m) + m (F - c(m)) = m F - (m-1) c(m)$ \\
\end{center}

\noindent which rearranged yields

\begin{center}
    $F/E \leq 2/m + (1-1/m) c(m)/E$.
\end{center}

\noindent We can now use the result of lemma \ref{thm:one} to bound the number of cycles of length $m$ or less, which is also a bound on $c(m)$ due to the fact that every face is a cycle. This gives

\begin{center}
    $F/E \leq 2/m + (1 - 1/m) c(m) / E \leq 2/m + o(1)$.
\end{center}

This holds for any constant $m$ meaning $F/E = o(1)$ w.h.p. Consequently, in combination with \eqref{eq:one}, we have

\begin{center}
    $\alpha \geq \frac{d - 2}{4} - o(1)$
\end{center}

\noindent hence the proof. \hfill $\square$

\section{Numerical Experiments}
\begin{figure}[htb]
    \centering
    \includegraphics[width=0.4\textwidth]{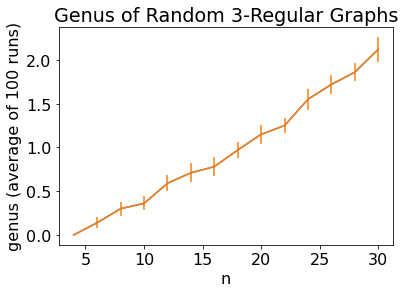}
    \caption{The growth of the genus of random 3-regular graphs}
    \label{fig:gen}
\end{figure}

Numerical experiments show that for small $n$, the genus grows much slower than $\frac{(d - 2)n}{4}$. Figure \ref{fig:gen} shows that for random 3-regular graphs, $\alpha$ is around $\frac{1}{15}$, significantly less than the expected asymptotic $\frac{1}{4}$.

The data in figure \ref{fig:gen} was calculated by a rotational embedding scheme that was implemented in the functional programming language Haskell. The details of this algorithm are described in \cite{boothby2007}. A rotation system describes a cyclic ordering of neighbors of each vertex that describes an embedding of the graph. The algorithm for determining genus via a rotational embedding scheme involves brute-forcing every possible cyclic embedding around each node, resulting in different permutations of orderings around each vertex that are homeomorphic to the sphere, also called 2-cells (meaning that they themselves are planar), which can then be used to find minimum genus \cite{boothby2007}. Due to the NP-hard nature of this algorithm, it can only be practically run on small graphs. In this case graphs up to $n=30$ could be generated, as seen below.

\section*{Acknowledgements}
I am highly grateful to Cristopher Moore of the Santa Fe Institute for his mentorship.

\bibliographystyle{amsplain}
\bibliography{bibliography.bib}

\end{document}